\documentclass[11pt,letterpaper]{amsart}
\usepackage{amssymb,mathtools,microtype}
\usepackage{xcolor}
\definecolor{refblue}{RGB}{24,63,130}
\usepackage[colorlinks=true,linkcolor=refblue,citecolor=refblue,urlcolor=refblue]{hyperref}
\usepackage{placeins}
\usepackage{aliascnt}
\usepackage[nameinlink,noabbrev]{cleveref}
\hypersetup{pdftitle={An extremal unimodular lattice of rank 37},
 pdfauthor={Eric Kubischta}}
 
\newtheorem{theorem}{Theorem}[section]
\newtheorem*{maintheorem}{Main theorem}

\newaliascnt{proposition}{theorem}
\newtheorem{proposition}[proposition]{Proposition}
\aliascntresetthe{proposition}

\newaliascnt{lemma}{theorem}
\newtheorem{lemma}[lemma]{Lemma}
\aliascntresetthe{lemma}

\newaliascnt{corollary}{theorem}
\newtheorem{corollary}[corollary]{Corollary}
\aliascntresetthe{corollary}

\theoremstyle{remark}
\newaliascnt{remark}{theorem}

\aliascntresetthe{remark}

\crefname{theorem}{Theorem}{Theorems}
\crefname{proposition}{Proposition}{Propositions}
\crefname{lemma}{Lemma}{Lemmas}
\crefname{corollary}{Corollary}{Corollaries}
\crefname{section}{Section}{Sections}

\title[An extremal unimodular lattice of rank 37]
{An extremal unimodular lattice of rank 37}
\author{Eric Kubischta}
\address{Department of Mathematics, Florida State University,
Tallahassee, FL 32306, USA}
\email{emk25g@fsu.edu}

\subjclass[2020]{Primary 11H06; Secondary 11H50, 11H56, 94B05}
\keywords{extremal unimodular lattice, minimum norm,
Cayley table, discriminant form, self-dual code}
\date{\today}

\newcommand{\Z}{\mathbb Z}
\newcommand{\R}{\mathbb R}
\newcommand{\F}{\mathbb F}
\DeclareMathOperator{\Aut}{Aut}
\DeclareMathOperator{\AGL}{AGL}
\numberwithin{equation}{section}
\allowdisplaybreaks[1]

\begin{document}

\begin{abstract}
We construct an extremal unimodular lattice of rank $37$,
resolving the last undetermined optimal minimum in ranks
at most $48$. The construction begins with the Cayley table
of $\Z/13\Z$, assigning a vector of norm four to each cell.
We prove that the resulting incidence lattice has exactly
$26$ unimodular completions of minimum four, all isometric.
The resulting Cayley lattice $\Lambda$ has automorphism group
$\AGL(1,13)\times C_6$ and is $s$-extremal, with shadow
minimum $13/4$ and kissing number $35\,594$.
An explicit $4$-frame in $\Lambda$ yields a Type~I self-dual
$\Z/4\Z$-code of length $37$ and minimum Euclidean weight
$16$, resolving the corresponding open case in coding theory.
\end{abstract}

\maketitle


\section{Introduction}\label{sec:introduction}

Determining the largest possible minimum norm of a
positive-definite integral unimodular lattice in each rank
is a classical problem connecting quadratic forms, sphere
packing, modular forms, and coding theory.
Since a unimodular lattice has covolume one, maximizing its
minimum is equivalent to maximizing its sphere-packing density
within this arithmetic class. The root lattice $E_8$ and the
Leech lattice are familiar examples of optimal solutions
\cite{CS}.

Modular forms provide general upper bounds, but attaining them
is a separate existence problem. The bound of Mallows, Odlyzko
and Sloane for even unimodular lattices \cite{MOS}, extended to
odd lattices by Rains and Sloane using shadow theory
\cite[Theorem~1]{RS}, gives
\[
\min K\leq 2\left\lfloor\frac{n}{24}\right\rfloor+2
\qquad(n\ne23)
\]
for an integral unimodular lattice $K$ of rank $n$, with
exceptional bound three in rank $23$.
A lattice attaining this bound is called \emph{extremal},
while one attaining the largest minimum actually occurring
in its rank is called \emph{optimal}. These notions need
not coincide: an upper bound may be unattainable, and
determining the optimal minimum requires both constructions
and nonexistence theorems.

Conway and Sloane determined the optimal minimum in every
rank through $33$, announcing the results in
\cite[Theorem~2]{CS1990} and giving a detailed treatment
in \cite{CSoptimal}. In particular, they proved that
minimum four is impossible in rank $33$, despite being
allowed by the general bound. Their 1998 table still
recorded several subsequent ranks, including $37$, as
undetermined between three and four
\cite[Table~1]{CSoptimal}.
Bachoc, Nebe and Venkov then excluded minimum four in
ranks $34$ and $35$ \cite[Theorem~1]{BNV}.
Since minimum-four lattices were already known in ranks
$36$, $38$, $39$, and $40$, they explicitly identified
rank $37$ as the only remaining case through rank $40$.

Further constructions, including Gaborit's code-based
constructions \cite{Gaborit2004} and Harada's work in
rank $41$ \cite{HaradaShadow,HaradaZ4}, settled the
surrounding higher ranks. Rank $37$ was consequently
left as the sole undetermined optimal minimum among
integral unimodular lattices of rank at most $48$.
The Nebe--Sloane catalogue continued to record its
optimal minimum as three or four \cite{NSCatalogue}.
Thus the unresolved question was whether rank $37$
admits an extremal lattice, or instead requires a
further nonexistence theorem like those in ranks
$33$, $34$, and $35$.
Because every smaller odd rank has optimal minimum
at most three, an affirmative answer would also
identify the smallest odd rank admitting a unimodular
lattice of minimum four.

The same gap appeared in coding theory. Construction A
over $\Z/4\Z$ converts a Type~I self-dual code of
minimum Euclidean weight $16$ into an odd unimodular
lattice of minimum four. Harada determined the largest
minimum Euclidean weight of Type~I self-dual
$\Z/4\Z$-codes in every length through $47$ except
$37$ \cite{HaradaZ4}. At length $37$, codes of minimum
Euclidean weight $12$ were known, but extensive searches
did not produce one of minimum Euclidean weight $16$
\cite{HaradaZ4}.
The lattice and code questions are closely related,
but the existence of a minimum-four lattice alone does
not supply such a code: a $4$-frame is also required.

We resolve both existence questions.

\begin{maintheorem}
There exists an integral unimodular lattice $\Lambda$
of rank $37$ and minimum four, with
\[
\Aut(\Lambda)\cong\AGL(1,13)\times C_6.
\]
It is $s$-extremal, with shadow minimum $13/4$ and
kissing number $35\,594$, and contains a $4$-frame.
Consequently, there exists a Type~I self-dual
$\Z/4\Z$-code of length $37$ and minimum Euclidean
weight $16$.
\end{maintheorem}

The key to our construction is that the Cayley table of
$\Z/13\Z$ reduces a large short-vector problem to a single
finite-geometric obstruction. Assigning norm-four vectors
to its cells reveals the root lattice $A_{12}^3$, whose
weight-lattice geometry shows that the only vectors of
norm below four in the relevant dual lattice are roots.
Remarkably, the Cayley relations force all these root
classes onto a single isotropic line
(\Cref{lem:short-dual,lem:root-line}).
Thus extremality reduces to choosing a unimodular glue
plane avoiding this line. The table symmetries make
all such choices isometric
(\Cref{thm:gluings,prop:gluing-isometry}),
giving the \emph{Cayley lattice} $\Lambda$.

The symmetry also singles out a characteristic vector
of norm thirteen by averaging the cell vectors.
This is precisely the norm needed to attain Gaborit's
shadow bound, and excluding shorter characteristic
vectors establishes $s$-extremality (\Cref{prop:shadow}).
More subtly, the same vector allows us to recover the
original Cayley configuration intrinsically from $\Lambda$.
Thus the symmetries visible in the construction,
together with negation, account for the full automorphism
group (\Cref{thm:automorphisms}).

Finally, the connection with coding theory emerges from
the completed lattice itself. An explicit $4$-frame
reveals a Type~I self-dual $\Z/4\Z$-code
(\Cref{prop:frame,cor:code}), whose coordinates are
given in \Cref{app:frame}.

\subsection*{Computational reproducibility}
An explicit Gram matrix and Magma verification code for the
Cayley lattice are available at
\url{https://github.com/Eric-Kubischta/extremal-unimodular-37}.
The code reconstructs $\Lambda$ from its defining vectors in
$\R^{40}$ and verifies its main properties, including extremality,
shadow, $4$-frame, associated code, and automorphism group.
A shorter verification starts directly from the Gram matrix.

\subsection*{Declaration of AI Use}
The author used ChatGPT (OpenAI), including GPT-6 and ChatGPT Work,
as an interactive aid in exploratory research and manuscript preparation.
This included computational searches that identified candidate lattices,
exploration of combinatorial constructions and intermediate arguments,
refinement of the exposition, and assistance in developing the
accompanying Magma verification code. AI-assisted suggestions were
treated as provisional and critically evaluated. The author
substantially developed, revised, and independently verified the
final mathematical arguments and takes full responsibility for the
results, computations, references, and text.

\section{Cayley Construction}\label{sec:incidence}

Let $e$, together with $e_i^R,e_i^C,e_i^S$ for $i\in\Z/13\Z$,
be an orthonormal basis of $\R^{40}$. All indices are taken modulo
$13$. The three blocks correspond to the \underline{R}ows,
\underline{C}olumns and \underline{S}ymbols of the Cayley table of
$\Z/13\Z$. The cell $(i,j)$ contains $i+j$, so its incidence
vector is $e_i^R+e_j^C+e_{i+j}^S$. To obtain norm four, we
augment these norm-three vectors by $e$ and set
\[
v_{ij}=e+e_i^R+e_j^C+e_{i+j}^S.
\]
Their inner products record the table incidences:
\begin{equation}\label{eq:ip}
(v_{ij},v_{k\ell})=1+\delta_{ik}+\delta_{j\ell}+\delta_{i+j,k+\ell}.
\end{equation}
Two distinct cells share at most one row, column or symbol, so
these inner products are $4$ on the diagonal, $2$ for incident
cells and $1$ otherwise. Let $V$ and $L$ be the real and integer
spans of the $v_{ij}$, respectively.

\begin{lemma}\label{lem:rank}
$L$ is an even integral lattice of rank $37$ in $V$.
\end{lemma}
\begin{proof}
Averaging over the cells in row $i$ gives
\[
\frac1{13}\sum_j v_{ij}
=e_i^R+\left(e+\frac1{13}\sum_j(e_j^C+e_j^S)\right).
\]
The expression in parentheses is independent of $i$, so
subtracting two row averages gives $\frac1{13}\sum_j(v_{ij}-v_{kj})=e_i^R-e_k^R$. The same argument applies to columns and to the cells carrying
a fixed symbol. Thus, for $X=R,C,S$,
\[
e_i^X-e_k^X\in V,
\qquad
13(e_i^X-e_k^X)\in L.
\]

It follows that every vector in $V^\perp$ has constant
coordinates within each block. Such a vector is orthogonal
to all $v_{ij}$ precisely when its three block values and
its $e$-coordinate sum to zero. These four values satisfy
one independent linear relation, so $\dim V^\perp=3$ and $\dim V=37$. Every $v_{ij}$ has coordinate sum in each block equal to its
$e$-coordinate. These three independent linear conditions
define a $37$-dimensional subspace containing $V$, so they
describe $V$ exactly.

Finally, $L\subseteq\Z^{40}$ is discrete and spans $V$, so
it has rank $37$. By \eqref{eq:ip}, $L$ is integral, and it
is even because its generators have norm four.
\end{proof}

An even lattice of odd rank cannot be unimodular, so we seek an
odd integral enlargement. Let $G$ be the full table symmetry group,
allowing the row, column and symbol roles to be interchanged.
These symmetries permute the $v_{ij}$ and preserve their inner
products, so $G\leq\Aut(L)$.

To retain these symmetries, we first look for a vector fixed by
$G$. Since $G$ acts transitively on the spanning configuration,
averaging shows that its fixed space in $V$ is the line through
\[
\bar v=\frac1{13^2}\sum_{i,j}v_{ij}
=e+\frac1{13}\sum_i(e_i^R+e_i^C+e_i^S).
\]
Now $(\bar v,v_{ij})=\bar v^2=16/13$, so adjoining $t\bar v$
preserves integrality precisely when $16t/13$ and $16t^2/13$ are
integers. Writing $t=13m/16$, the second condition becomes
$4\mid m$. The shortest nonzero choice is therefore
\[
\kappa =\frac{13}{4}\bar v
=\frac14\left(13e+\sum_i(e_i^R+e_i^C+e_i^S)\right).
\]
\begin{lemma}\label{lem:gamma}
$\Gamma=L+\Z \kappa$ is an odd integral overlattice of $L$ of index four,
preserved by $G$, and $\kappa$ is characteristic in $\Gamma$.
\end{lemma}
\begin{proof}
We have $(\kappa,v_{ij})=4$ and $\kappa^2=13$. The diagonal cells contain each row, column and symbol once,
so $4\kappa=\sum_i v_{ii}\in L$. On the other hand, $m \kappa\notin L$
for $m=1,2,3$, since its block coordinates are $m/4$, whereas
$L\subseteq\Z^{40}$. Hence $[\Gamma:L]=4$.
The group $G$ preserves $L$ and fixes $\kappa$. Finally,
$(\kappa,x)\equiv x^2\pmod2$ holds on the generators $v_{ij},\kappa$,
and hence on $\Gamma$.
\end{proof}

However, $\Gamma$ is not unimodular, as the following lemma shows.

\begin{lemma}\label{lem:det}
We have $\det L=16\cdot13^4$ and $\det\Gamma=13^4$.
\end{lemma}
\begin{proof}
Put $P=V\cap\Z^{40}$. By \Cref{lem:rank}, the coordinate
differences lie in $P$. Moreover, subtracting $(x,e)v_{00}$
from any $x\in P$ leaves an integral vector whose coordinate
sum in each block is zero. Thus
\[
P=\Z v_{00}\oplus A,\qquad
A=A_{12}^R\oplus A_{12}^C\oplus A_{12}^S,
\]
where $A_{12}^X$ is the root lattice of integral vectors
whose coordinates in block $X$ sum to zero.
Since $\det A=13^3$ and the component of $v_{00}$
orthogonal to $A$ is $\bar v$, we obtain $\det P=13^3\bar v^2=16\cdot13^2$. 

To compute $[P:L]$, write
$r_i=e_i^R-e_0^R$, $c_i=e_i^C-e_0^C$, and
$s_i=e_i^S-e_0^S$ so that
\[
v_{ij}=v_{00}+r_i+c_j+s_{i+j}.
\]
Then $v_{i0} - v_{00} = r_i + s_i \in L$ and $v_{0i} - v_{00} = c_i + s_i \in L$. Thus
\[
    r_i \equiv c_i \equiv -s_i \pmod{L}. 
\]
Moreover, $v_{ij}-v_{i0}-v_{0j}+v_{00}
=s_{i+j}-s_i-s_j\in L$ so
\[
s_{i+j}\equiv s_i+s_j\pmod L.
\]
Since $P$ is generated by $v_{00}$ and the $r_i,c_i,s_i$,
while $L$ is generated by the $v_{ij}$, these relations
give a presentation of $P/L$. They reduce it to a single
generator $s_1$, with the relation $13s_1=0$.
Hence $P/L\cong\Z/13\Z$.

Therefore $[P:L]=13$, and the index--determinant formula gives both $\det L=13^2\det P=16\cdot13^4$ and $\det\Gamma=\frac{\det L}{4^2}=13^4$. 
\end{proof}

Symmetry motivated adjoining $\kappa $ as a first step toward a unimodular lattice. Although the resulting $\Gamma$ is not yet unimodular, this step is in fact unavoidable.

\begin{lemma}\label{lem:canonical}
Every integral unimodular overlattice of $L$ contains $\Gamma$.
\end{lemma}
\begin{proof}
Since $(\kappa/4,v_{ij})=1$, the vector $\kappa/4$ lies in the dual lattice $L^*$.
Its class has order $16$: its block coordinates are $1/16$,
and $4\kappa \in L$. Since $|L^*/L|=16\cdot13^4$, its $2$-primary
part is cyclic of order $16$. Its unique
subgroup of order four is $\Gamma/L$. For an integral unimodular overlattice $M$, we have
$M/L\subseteq L^*/L$ and
$|M/L|=\sqrt{\det L}=4\cdot13^2$.
Thus its $2$-primary part must be $\Gamma/L$,
and hence $\Gamma\subseteq M$.
\end{proof}
All remaining glue lies in the discriminant group
$D=\Gamma^*/\Gamma$. For $x\in\Gamma^*$, write
$[x]=x+\Gamma\in D$.

\begin{lemma}\label{lem:discriminant}
For every $x\in V$, we have the frame identity
\begin{equation}\label{eq:frame-operator}
\sum_{i,j}(x,v_{ij})v_{ij}=13x+15(x,\kappa)\kappa.
\end{equation}
Consequently, $D$ is a four-dimensional vector space over $\F_{13}$.
\end{lemma}
\begin{proof}
The identity follows from the averaging formulas in
\Cref{lem:rank}: both sides act by $13$ on the coordinate differences
within each block and send $\kappa$ to $208\kappa$; these vectors span $V$.
For $x\in\Gamma^*$, the sum lies in $L$ and $(x,\kappa)\in\Z$, so
$13x\in\Gamma$. Thus $13D=0$, making $D$ an $\F_{13}$-vector space,
and $|D|=\det\Gamma=13^4$ gives its dimension.
\end{proof}

Glue vectors must pair integrally with one another. Their
pairings are recorded by the nondegenerate discriminant form
\[
\beta:D\times D\longrightarrow\F_{13},\qquad
\beta([x],[y])=13(x,y)\pmod{13}.
\]
The inclusion $13\Gamma^*\subseteq\Gamma$ makes this well-defined,
and $\beta([x],[y])=0$ precisely when $(x,y)\in\Z$.
Nondegeneracy follows from $(\Gamma^*)^*=\Gamma$.

\begin{proposition}\label{prop:isotropic}
The integral unimodular overlattices of $L$ correspond bijectively
to the totally isotropic planes of $(D,\beta)$.
\end{proposition}
\begin{proof}
For $H\leq D$, its inverse image
\[
\Lambda_H=\{x\in\Gamma^*:[x]\in H\}
\]
is integral exactly when $\beta(H,H)=0$.
If $\dim H=r$, then
$[\Lambda_H:\Gamma]=13^r$ and $\det\Lambda_H=13^{4-2r}$,
so unimodularity is equivalent to $r=2$.
Every integral unimodular overlattice of $L$ contains
$\Gamma$ by \Cref{lem:canonical}.
\end{proof}

We call $H$ the \emph{glue plane}; adjoining representatives of
its basis gives $\Lambda_H$. The vector $\kappa$ remains characteristic
after every integral gluing: for $x\in\Lambda_H$, apply its
characteristic congruence in $\Gamma$ to $13x$ and reduce modulo
two to obtain $(\kappa,x)\equiv x^2\pmod2$.

\section{Extremal Gluing}

To obtain minimum four, we must exclude glue that introduces
vectors of norm one, two or three. Since every unimodular
gluing lies in $\Gamma^*$, we begin by determining its
vectors of these norms.

\begin{lemma}\label{lem:short-dual}
$\Gamma^*$ has no vectors of norm one or three, and its
norm-two vectors are precisely the roots of $A$, given by
$e_i^X-e_j^X$ for $X=R,C,S$ and $i\ne j$.
\end{lemma}
\begin{proof}
\emph{1. Reduction to $N^*$.}
Let $x\in\Gamma^*$ have $x^2\in\{1,2,3\}$.
Since $V=\R \kappa \perp\operatorname{span}_{\R}A$ and
$(x,\kappa )\in\Z$, we may write
\[
    x=\frac{m}{13}\kappa +u
\]
where $m=(x,\kappa )\in\Z$ and $u\in\operatorname{span}_{\R}A$. It follows that
\begin{equation}
    13x^2= m^2 + 13u^2. \label{eqn:shortvec}
\end{equation}
Put $N=L\cap A$. Since $N\subseteq\Gamma$ and $\kappa \perp A$,
we have $(u,N)=(x,N)\subseteq\Z$, so $u\in N^*$,
with the dual taken in $\operatorname{span}_{\R}A$.

Choose the positive roots to be $e_j^X-e_i^X$ for $0\le i<j\le12$
in the row and column blocks, and their negatives in the
symbol block. The corresponding Weyl vector of $A$, half the sum
of these positive roots, is
\[
\rho=\sum_{i=0}^{12}(i-6)(e_i^R+e_i^C-e_i^S).
\]
The signs encode the Cayley relation $i+j-s=0$.
Indeed, its pairings with $r_i,c_i,s_i$ are $i,i,-i$,
so the quotient relations in \Cref{lem:det} give
\[
N=\{a\in A:(a,\rho)\equiv0\pmod{13}\},
\qquad [A:N]=13.
\]
Thus $\rho/13\in N^*$, and $(\rho/13,r_1)=1/13$
shows that its class generates $N^*/A^*$. Consequently,
\[
N^*=A^*+\Z\frac{\rho}{13},
\qquad
u=w+\frac{k}{13}\rho
\quad(w\in A^*,\ k\in\Z).
\]

\smallskip
\noindent
\emph{2. Deep holes force $u\in A^*$.}
Each block of $\rho/13$ is a deep hole of $A_{12}^*$,
at squared distance $14/13$ from that lattice
\cite[Chapter~21, \S3F]{CS}.
If $13\nmid k$, multiplication by $k$ permutes the
centered residues modulo $13$, so $k\rho/13$ differs
from a blockwise coordinate permutation of $\rho/13$
by an element of $A$.
Squared distances add across the three orthogonal blocks,
giving
\begin{equation}\label{eq:coset-bound}
u^2\ge
\operatorname{dist}\left(\frac{k\rho}{13},A^*\right)^2
=3\cdot\frac{14}{13}
=\frac{42}{13}>3.
\end{equation}
This contradicts $u^2\le x^2\le3$. Hence $u\in A^*$.

\smallskip
\noindent
\emph{3. Weight shells leave only roots.}
For each block $X$, let $a_X\in\{0,\ldots,12\}$ label
the class of $u_X$ in $A_{12}^*/A_{12}\cong\Z/13\Z$,
with coordinates congruent to $-a_X/13$ modulo $\Z$.
The minimum norm in class $a_X$ is
$a_X(13-a_X)/13$ \cite[Chapter~4, \S6]{CS}.
Since norms within a coset differ by even integers,
\[
13u_X^2=a_X(13-a_X)+26k_X,
\qquad k_X\in\Z_{\ge0}.
\]
Summing the three block contributions and retaining
only $u^2\le3$, the possible values are
\[
13u^2\in\{0,12,22,24,26,30,34,36,38\}.
\]
Put $a=a_R+a_C+a_S\pmod{13}$.
Comparing these with the condition in \cref{eqn:shortvec} leaves only five shells:
\[
\begin{array}{c|c|c|c}
x^2&|m|&13u^2&a\\ \hline
1&1&12&\pm1\\
2&0&26&0\\
2&2&22&\pm2\\
3&1&38&\pm1\\
3&3&30&\pm3
\end{array}
\]
In particular, $a\equiv\pm m\pmod{13}$. Since $u_0^X\equiv-a_X/13\pmod{\Z}$ and
$(\kappa ,v_{00})=4$, integrality of $(x,v_{00})$ gives
\[
(x,v_{00})\equiv\frac{4m-a}{13}\pmod{\Z},
\qquad\text{hence}\qquad a\equiv4m\pmod{13}.
\]
So $4m\equiv\pm m\pmod{13}$, impossible
for $0<|m|\le3$. Therefore $m=0$, and the remaining row gives
$x=u\in A$ with $x^2=2$. These are precisely the stated
roots. Conversely, every such root pairs integrally
with the $v_{ij}$ and is orthogonal to $\kappa $, so belongs
to $\Gamma^*$.
\end{proof}

The only obstruction to minimum four is therefore a root.
All $468$ roots are controlled by one line of glue classes.

\begin{lemma}\label{lem:root-line}
The root classes in $D$ are precisely the nonzero elements of
the isotropic line $\ell=\langle[s_1]\rangle$.
\end{lemma}
\begin{proof}
By \Cref{lem:short-dual}, every root is of the form
$r_i-r_j$, $c_i-c_j$, or $s_i-s_j$ for $i\ne j$.
The relations in \Cref{lem:det}, together with
$L\subseteq\Gamma$, give
\begin{equation}\label{eq:root-classes}
[r_i]=[c_i]=-[s_i]=-i[s_1].
\end{equation}

We have $\Gamma\cap P=L$, since a multiple of $\kappa$
has integral coordinates only if it belongs to $L$.
As $s_1$ generates $P/L\cong\Z/13\Z$ by \Cref{lem:det},
we have $[s_1]\ne0$.

Thus every root represents a nonzero multiple of $[s_1]$,
and the roots $s_i$ for $i\ne0$ represent all such multiples.
Finally,
\[
\beta([s_1],[s_1])=13s_1^2=26\equiv0\pmod{13},
\]
so $\ell$ is isotropic.
\end{proof}

A glue plane $H$ gives minimum four exactly when
$H\cap\ell=0$. To find these planes, consider
$\Sigma=\ell^\perp/\ell$ with its induced form $\beta_\Sigma$.
Since $D$ is nondegenerate of dimension four and $\ell$ is
isotropic, $\ell^\perp$ has dimension three with radical $\ell$.
Thus $\Sigma$ is a nondegenerate two-dimensional space.

\begin{lemma}\label{lem:sigma}
The space $(\Sigma,\beta_\Sigma)$ is a hyperbolic plane.
\end{lemma}
\begin{proof}
For $X=R,C,S$, consider the natural block weights
\[
\omega_X=e_0^X-\frac1{13}\sum_i e_i^X,
\]
which are minuscule weights of $A_{12}$ (the first fundamental
weights for the usual ordering of the coordinates). They are mutually orthogonal
of squared norm $12/13$.

For integers $a_X$, we have
\[
\begin{aligned}
\left(\sum_X a_X\omega_X,v_{ij}\right)
&=a_R\delta_{i0}+a_C\delta_{j0}
+a_S\delta_{i+j,0} -\frac1{13}\sum_X a_X.
\end{aligned}
\]
Since each $\omega_X$ is orthogonal to $\kappa$, the
combination belongs to $\Gamma^*$ precisely when $\sum_X a_X\equiv0\pmod{13}$. Its pairing with $s_1$ is $-a_S$, so its class also lies
in $\ell^\perp$.

For two such combinations, mutual orthogonality and
$\omega_X^2=12/13$ give
\[
\begin{aligned}
\beta\left(
\left[\sum_X a_X\omega_X\right],
\left[\sum_X b_X\omega_X\right]
\right)
&\equiv12\sum_X a_Xb_X \equiv-\sum_X a_Xb_X\pmod{13}.
\end{aligned}
\]
Thus the discriminant pairing on these coefficient triples
is the negative of the ordinary dot product over $\F_{13}$.
In particular, such a class is isotropic precisely when $\sum_X a_X^2\equiv0\pmod{13}$.

We therefore seek nonzero solutions of $\sum_X a_X=\sum_X a_X^2=0$ in $\F_{13}$. Necessarily $a_R\ne0$, since otherwise
$a_C=-a_S$ and $2a_C^2=0$.
After rescaling, we may take
$(a_R,a_C,a_S)=(1,t,-1-t)$, giving
\[
1+t^2+(1+t)^2=0,
\qquad\text{or}\qquad
t^2+t+1=0.
\]
The two solutions are $t=3,9$. Convenient integer
representatives of these two directions are
$(1,3,-4)$ and $(3,1,-4)$, so we take
\[
\theta=\omega_R+3\omega_C-4\omega_S,
\qquad
\theta'=3\omega_R+\omega_C-4\omega_S.
\]
Notice that $\theta'$ is obtained from $\theta$ by
interchanging the row and column blocks.

Since $(\theta,s_1)=(\theta',s_1)=4$,
both classes lie in $\ell^\perp$.
Their coefficient triples have squared length $26$
and mutual dot product $22$, so their images in
$\Sigma$ have Gram matrix $\begin{psmallmatrix}
0&4\\
4&0
\end{psmallmatrix}$ with respect to $\beta_\Sigma$.

This matrix is nondegenerate, so the two images form
a basis of $\Sigma$, exhibiting it as a hyperbolic plane.
Its only isotropic lines are the two basis lines,
since a vector with coordinates $u,v$ has self-pairing
$8uv$.
\end{proof}

\begin{theorem}\label{thm:gluings}
$\Gamma$ has exactly $28$ integral unimodular gluings:
two of minimum two and $26$ of minimum four.
\end{theorem}
\begin{proof}
The hyperbolic plane spanned by $[\theta],[\theta']$
has a nondegenerate orthogonal complement containing
the isotropic vector $[s_1]$.
We may therefore choose $\tau\in\Gamma^*$ whose class
lies in this complement and satisfies $\beta([s_1],[\tau])=4$. 

Adding a suitable multiple of $s_1$ makes $[\tau]$
isotropic, since this changes its self-pairing by
a multiple of $8$. Thus
\[
D=\langle[s_1],[\tau]\rangle
\perp\langle[\theta],[\theta']\rangle,
\]
where the only nonzero pairings between basis vectors are $\beta([s_1],[\tau])
=\beta([\theta],[\theta'])=4$. 

A totally isotropic plane $H$ containing $\ell$
descends to an isotropic line in
$\Sigma=\ell^\perp/\ell$.
By \Cref{lem:sigma}, there are exactly two such planes: $\langle[s_1],[\theta]\rangle$ and $\langle[s_1],[\theta']\rangle$. 

Now suppose $H\cap\ell=0$.
Since $H$ is a totally isotropic plane in a
nondegenerate four-dimensional space, $H=H^\perp$.
Thus $H$ cannot lie in $\ell^\perp$, as this would
imply $\ell\subseteq H$.
Consequently, $H\cap\ell^\perp$ is a line whose image
in $\Sigma$ is isotropic.
By \Cref{lem:sigma}, it has a generator of the form
\[
x=[\theta]+t[s_1]
\quad\text{or}\quad
x=[\theta']+t[s_1],
\qquad t\in\F_{13}.
\]

Consider the first case.
Choose a second basis vector $y\in H$ with
$[\tau]$-coefficient one.
Subtracting a multiple of $x$, we may write
\[
y=[\tau]+a[s_1]+d[\theta'].
\]
The conditions that $H$ be totally isotropic give
\[
0=\beta(x,y)=4(t+d),
\qquad
0=\beta(y,y)=8a.
\]
Hence $d=-t$ and $a=0$.
The second case is identical with $\theta$ and
$\theta'$ interchanged. Consequently, the remaining
planes are precisely
\begin{equation}\label{eq:glue-planes}
\begin{aligned}
H_t&=\langle[\theta]+t[s_1],\,
             [\tau]-t[\theta']\rangle,\\
H'_t&=\langle[\theta']+t[s_1],\,
              [\tau]-t[\theta]\rangle,
\end{aligned}
\qquad t\in\F_{13}.
\end{equation}
Conversely, these $26$ planes are totally isotropic
and avoid $\ell$. They are distinct because their
intersections with $\ell^\perp$ are the distinct lines
\[
\langle[\theta]+t[s_1]\rangle,
\qquad
\langle[\theta']+t[s_1]\rangle.
\]

Together with the two planes containing $\ell$,
this gives exactly $28$ integral unimodular gluings
by \Cref{prop:isotropic}.
Finally, \Cref{lem:short-dual,lem:root-line} show that
the two planes containing $\ell$ have minimum two,
while the remaining $26$ have minimum four.
\end{proof}

\begin{proposition}\label{prop:gluing-isometry}
The unimodular gluings of $\Gamma$ form exactly two
isometry classes, of minimum two and four.
\end{proposition}
\begin{proof}
Consider the Cayley-table translation
\[
U:(i,j)\longmapsto(i+1,j+1),
\]
which acts on the symbol indices by $s\mapsto s+2$.
It permutes the $v_{ij}$ and fixes $\kappa$, so
it preserves $\Gamma$.

By \eqref{eq:root-classes}, we have
\[
[Us_1]=[s_3-s_2]=[s_1].
\]
Moreover,
\[
\begin{aligned}
[U\theta]-[\theta]
&=[r_1+3c_1-4s_2] =(-1-3-8)[s_1]=[s_1],
\end{aligned}
\]
and similarly $[U\theta']=[\theta']+[s_1]$.

By \Cref{thm:gluings}, a glue plane avoiding $\ell$
is uniquely determined by its intersection with
$\ell^\perp$. Consequently,
\[
U(H_t)=H_{t+1},
\qquad
U(H'_t)=H'_{t+1}.
\]
Thus $U$ acts transitively on each family of
$13$ minimum-four glue planes.

Interchanging the row and column blocks also preserves
$\Gamma$, fixes $[s_1]$, and exchanges $[\theta]$ and
$[\theta']$. It therefore exchanges the two families,
as well as the two planes containing $\ell$.

Hence all $26$ minimum-four gluings are mutually
isometric, as are the two minimum-two gluings.
Since their minima differ, these give exactly
two isometry classes.
\end{proof}

We henceforth fix the minimum-four glue plane
\[
H_0=\langle[\theta],[\tau]\rangle
\]
and call $\Lambda=\Lambda_{H_0}$ the \emph{Cayley lattice}.
We now construct the second glue vector $\tau$ naturally
from the Weyl vector and fundamental weights.

\begin{lemma}\label{lem:explicit-glue}
We may take $\tau=-\frac4{13}\rho
-2(\omega_R+\omega_C-\omega_S)-s_1$. In particular,
\[
\Lambda=\Gamma+\Z\theta+\Z\tau.
\]
\end{lemma}
\begin{proof}
We seek a vector $\tau\in\Gamma^*$ whose class pairs
with $[s_1]$ as $4$, is orthogonal to $[\theta]$ and
$[\theta']$, and is isotropic.

The Weyl vector provides a natural starting point.
Since $(\rho,s_1)=-1$, the vector $-4\rho/13$ satisfies $13\left(-\frac4{13}\rho,s_1\right)=4$, giving the desired discriminant pairing. However, it does not
belong to $\Gamma^*$ and is not orthogonal to
$\theta,\theta'$. We correct both problems using
the fundamental weights, writing
\[
x=-\frac4{13}\rho+\sum_X a_X\omega_X,
\qquad a_X\in\Z.
\]
By the definitions of $\rho$ and $\omega_X$,
\[
(\rho,v_{ij})\equiv-6\pmod{13},
\qquad
(\omega_X,v_{ij})\equiv-\frac1{13}\pmod{\Z}.
\]
Since $x\perp\kappa$, the condition $x\in\Gamma^*$
is therefore equivalent to $\sum_X a_X\equiv24\pmod{13}$. 

Next, we impose the stronger conditions
$(x,\theta)=(x,\theta')=0$. Using
$(\rho,\theta)=(\rho,\theta')=-48$ and
$\omega_X^2=12/13$, these become
\[
\begin{aligned}
a_R+3a_C-4a_S&=-16, \qquad 
3a_R+a_C-4a_S&=-16.
\end{aligned}
\]
Thus $a_R=a_C$ and $a_S=a_R+4$.
Together with the integrality condition, this gives
\[
(a_R,a_C,a_S)\equiv(-2,-2,2)\pmod{13}.
\]
Choosing these integer representatives, we obtain
\[
x=-\frac4{13}\rho
-2(\omega_R+\omega_C-\omega_S)\in\Gamma^*,
\]
with $(x,\theta)=(x,\theta')=0$ and
$\beta([x],[s_1])=4$.

It remains only to make $[x]$ isotropic. Direct
calculation gives $x^2=\frac{528}{13}$ and $13(x,s_1)=-22$ so
\[
\beta([x],[x])=8,
\qquad
\beta([x],[s_1])=4.
\]
Since $[s_1]$ is isotropic and orthogonal to
$[\theta],[\theta']$, adding a multiple of $[s_1]$
preserves all the conditions already established.
Moreover,
\[
\beta([x]+t[s_1],[x]+t[s_1])=8+8t.
\]
Taking $t=-1$ therefore gives the isotropic class
\[
[\tau]=[x]-[s_1].
\]
Thus $\tau=x-s_1$ satisfies all four required
pairing conditions and represents the class chosen
in \Cref{thm:gluings}. These conditions uniquely determine the class
$[\tau]$ chosen in \Cref{thm:gluings}: any two such
classes differ by a multiple of $[s_1]$, and
isotropy forces that multiple to vanish.
\end{proof}

\begin{proposition}\label{prop:shadow}
The Cayley lattice $\Lambda$ is $s$-extremal. Its shadow
has minimum $13/4$, and its kissing number is $35\,594$.
\end{proposition}
\begin{proof}
Recall that the shadow of an odd unimodular lattice consists
of the halves of its characteristic vectors: $S(\Lambda)=\frac12\operatorname{Char}(\Lambda)$. Thus we seek the smallest norm of a characteristic vector.
By van der Blij's congruence \cite{vanDerBlij1959}, every characteristic
vector $d$ satisfies
\[
d^2\equiv\operatorname{rank}(\Lambda)
=37\equiv5\pmod8.
\]
Since $\kappa$ is characteristic and $\kappa^2=13$,
the only possible smaller norm is five.

Suppose, then, that $d$ is characteristic with $d^2=5$.
Since $d$ and $\kappa$ are both characteristic,
the vectors $(\kappa\pm d)/2$ belong to $\Lambda$.
They are nonzero, and their squared norms sum to $\frac{\kappa^2+d^2}{2}=9$. As $\Lambda$ has minimum four, one must have norm four.
Replacing $d$ by $-d$ if necessary, we may therefore take
\[
x=\frac{\kappa-d}{2},
\qquad x^2=4,\qquad (x,\kappa)=6.
\]
Put $t_{ij}=(x,v_{ij})\in\Z$.
Using $\sum_{i,j}v_{ij}=52\kappa$ and the frame identity
\eqref{eq:frame-operator},
\[
\sum_{i,j}t_{ij}=312,
\qquad
\sum_{i,j}t_{ij}^2=592.
\]
But $(t-1)(t-2)\ge0$ for every integer $t$, whereas
\[
\sum_{i,j}(t_{ij}-1)(t_{ij}-2)
=592-3(312)+2(169)=-6,
\]
a contradiction.

Therefore the smallest characteristic norm is $13$,
attained by $\kappa$, and hence
\[
\min S(\Lambda)=\frac{13}{4}.
\]

Gaborit's shadow bound \cite{Gaborit} states that an odd
unimodular lattice of rank $n$ satisfies
\[
8\min(\Lambda)+4\min S(\Lambda)\le n+8,
\]
apart from an exceptional case in rank $23$.
Equality defines an $s$-extremal lattice.
In our case,
\[
8\min(\Lambda)+4\min S(\Lambda)
=8\cdot4+4\cdot\frac{13}{4}
=45=37+8.
\]
Thus $\Lambda$ is $s$-extremal.

Finally, the standard theta-series identities for odd
unimodular lattices \cite{RS} give
\[
\Theta_\Lambda(q)
=\sum_{j=0}^{4}a_j\theta_3(q)^{37-8j}\Delta_8(q)^j,
\]
where $\theta_2,\theta_3,\theta_4$ are the usual Jacobi theta
functions and
\[
\Delta_8(q)=q\prod_{m\ge1}
(1-q^{2m-1})^8(1-q^{4m})^8.
\]
The corresponding shadow series is
\[
\Theta_{S(\Lambda)}(q)
=\sum_{j=0}^{4}\frac{(-1)^ja_j}{16^j}
\theta_2(q)^{37-8j}\theta_4(q^2)^{8j}.
\]
The absence of vectors of norm one, two, or three,
together with the shadow minimum $13/4$, determines
\[
(a_0,a_1,a_2,a_3,a_4)
=(1,-74,1036,-1184,0).
\]
Consequently,
\[
\Theta_\Lambda(q)=1+35\,594q^4+O(q^5),
\]
so the kissing number is $35\,594$.
\end{proof}

\section{Automorphisms}\label{sec:automorphisms}

The full table symmetry group $G$ preserves $\Gamma$, but
a symmetry preserves $\Lambda$ only when it stabilizes the
glue plane $H_0$. By \Cref{thm:gluings}, a minimum-four
glue plane is determined by its intersection with
$\ell^\perp$. Thus, for $H_0$, it suffices to preserve
the line $\langle[\theta]\rangle$.

Consider a table translation
$(i,j)\mapsto(i+a,j+b)$. By \eqref{eq:root-classes},
its action on the distinguished line is
\[
\begin{aligned}
[\theta]
&\longmapsto[\theta]+[r_a+3c_b-4s_{a+b}]=[\theta]-(5a+7b)[s_1].
\end{aligned}
\]
Hence the translations preserving $H_0$ are exactly
those with $b=3a$. We take
\[
T:(i,j)\longmapsto(i+1,j+3).
\]
Multiplication by a nonzero scalar fixes the block
weights $\omega_X$, while cyclically permuting the
row, column and symbol roles preserves
$\langle[\theta]\rangle$. Accordingly, consider
\[
M:(i,j)\longmapsto(2i,2j),
\qquad
C:(i,j)\longmapsto(j,-i-j).
\]
Indeed, $M\theta=\theta$, while
\[
C\theta-3\theta=-13\omega_C+13\omega_S\in L,
\]
where the inclusion follows from the relations in
\Cref{lem:det}. Thus $T,M,C$ preserve $\Lambda$.
They all fix $\kappa$. Moreover, $T$ cyclically permutes
the coordinates within each block, so its fixed
space in $V$ is precisely $\R\kappa$.

To prove that these symmetries generate the full
automorphism group, we first recover $\kappa$
intrinsically.

\begin{lemma}\label{lem:characteristic-line}
Every automorphism of $\Lambda$ sends $\kappa$ to
$\kappa$ or $-\kappa$.
\end{lemma}
\begin{proof}
We first show that no $x\in\Lambda$ satisfies
$x^2=(x,\kappa)=8$. Put $t_{ij}=(x,v_{ij})\in\Z$.
Since $\sum_{i,j}v_{ij}=52\kappa$, the frame identity
\eqref{eq:frame-operator} gives
\[
\sum_{i,j}t_{ij}=416,
\qquad
\sum_{i,j}t_{ij}^2=1064.
\]
But $(t-2)(t-3)\ge0$ for every integer $t$, whereas
\[
\sum_{i,j}(t_{ij}-2)(t_{ij}-3)
=1064-5(416)+6(169)=-2,
\]
a contradiction.

Let $\Lambda_0$ be the even sublattice of $\Lambda$
and consider the characteristic-vector shell
\[
\mathcal X=\{d\in\kappa+2\Lambda_0:d^2=13\}.
\]
For $d=\kappa-2x\in\mathcal X$, we have
$x\in\Lambda_0$ and
\[
x^2=(\kappa,x).
\]
If $x\ne0$, the minimum-four condition applied to
$x$ and $\kappa-x$ gives
$4\le x^2\le9$. Since $x^2$ is even and the
value eight was excluded, we have
\[
x^2\in\{0,4,6\},
\qquad
(\kappa,d)=13-2x^2\in\{13,5,1\}.
\]
Thus every vector in $\mathcal X$ pairs positively
with $\kappa$.

The translation $T$ preserves $\mathcal X$, so
\[
z=\sum_{d\in\mathcal X}d
\]
is fixed by $T$ and therefore belongs to $\R\kappa$.
Since $\kappa\in\mathcal X$ and $(\kappa,z)>0$,
we have $z\ne0$.

Finally, the characteristic vectors of $\Lambda$
form the disjoint union
\[
\kappa+2\Lambda
=(\kappa+2\Lambda_0)\sqcup(-\kappa+2\Lambda_0).
\]
Every automorphism preserves or exchanges these two
classes, and hence sends $\mathcal X$ to
$\mathcal X$ or $-\mathcal X$. It therefore sends
$z$ to $\pm z$, and consequently
$\kappa$ to $\pm\kappa$.
\end{proof}

Having recovered $\kappa$, we can recognize the
original Cayley configuration among the norm-four
vectors.

\begin{lemma}\label{lem:recover-incidence}
The vectors $x\in\Lambda$ satisfying
$x^2=(x,\kappa)=4$ are precisely the $v_{ij}$.
\end{lemma}
\begin{proof}
Write
\[
x=\frac4{13}\kappa+u,
\qquad u\perp\kappa.
\]
As in \Cref{lem:short-dual}, we have $u\in N^*$,
and the norm condition gives $13u^2=36$.
Since $36<42$, the deep-hole bound
\eqref{eq:coset-bound} forces $u\in A^*$.

The weight-shell formula in \Cref{lem:short-dual}
shows that the only possibilities at $13u^2=36$
are a single block of charge $\pm4$, or three
blocks of charge $\pm1$. The integrality condition
from that lemma gives
\[
a_R+a_C+a_S\equiv4(x,\kappa)
=16\equiv3\pmod{13}.
\]
This leaves only
\[
a_R=a_C=a_S=1.
\]
Thus each block is a shortest weight of charge one,
and we obtain
\[
x=e+e_i^R+e_j^C+e_k^S
\]
for some $i,j,k\in\Z/13\Z$.

But
\[
x-v_{ij}=e_k^S-e_{i+j}^S\in\Lambda
\]
would have norm two unless $k=i+j$, contrary
to the minimum of $\Lambda$. Hence $x=v_{ij}$.
The converse follows from
$v_{ij}^2=(v_{ij},\kappa)=4$.
\end{proof}

\begin{theorem}\label{thm:automorphisms}
The full automorphism group of the Cayley lattice is
\[
\Aut(\Lambda)=\langle T,M,C,-I\rangle
\cong\operatorname{AGL}(1,13)\times C_6,
\]
of order $936$.
\end{theorem}
\begin{proof}
Let $g\in\Aut(\Lambda)$. By
\Cref{lem:characteristic-line}, composing with
negation if necessary, we may assume $g\kappa=\kappa$.
Then \Cref{lem:recover-incidence} shows that
$g$ permutes the $169$ vectors $v_{ij}$.
Consequently, it preserves $L$, $\Gamma$ and
$\Gamma^*$.

By \Cref{lem:short-dual}, the roots of $\Gamma^*$
form the three orthogonal irreducible root systems
$A_{12}^R,A_{12}^C,A_{12}^S$. Thus $g$ permutes
these components. The projections of the $v_{ij}$
onto each component are precisely the thirteen
weights
\[
e_i^X-\frac1{13}\sum_k e_k^X,
\qquad i\in\Z/13\Z.
\]
These form a regular simplex and recover the
coordinate labels in each block. Since $g$ preserves
the incidence vectors, it induces permutations
of the rows, columns and symbols preserving the
Cayley relation $i+j=s$. Therefore $g\in G$.

A table symmetry preserving the three roles
satisfies
\[
\gamma(i+j)=\alpha(i)+\delta(j).
\]
Subtracting the values at zero shows that the
three permutations have a common additive part.
Thus every such symmetry has the form
\[
(i,j)\longmapsto(ai+b,aj+d),
\qquad
a\in\F_{13}^{\times},\quad b,d\in\F_{13}.
\]
The six permutations of the signed roles
$(i,j,-i-j)$ give the remaining symmetries.
Hence
\[
|G|=13^2\cdot12\cdot6.
\]

By \Cref{prop:gluing-isometry}, $G$ acts transitively
on the $26$ minimum-four glue planes. Therefore
\[
|\Aut(\Lambda)_\kappa|
=|\operatorname{Stab}_G(H_0)|
=\frac{13^2\cdot12\cdot6}{26}
=468.
\]

The generators already constructed satisfy
\[
M T M^{-1}=T^2,\qquad
C T C^{-1}=T^3,\qquad
MC=CM.
\]
Thus $\langle T,M\rangle\cong\operatorname{AGL}(1,13)$.
Moreover, $M^8C$ is an independent central element
of order three, since $2^8\cdot3\equiv1\pmod{13}$.
Consequently,
\[
\langle T,M,C\rangle
\cong\operatorname{AGL}(1,13)\times C_3
\]
has order $468$ and is the full stabilizer of
$\kappa$. Adjoining the independent central
involution $-I$ proves the theorem.
\end{proof}

\section{A Self-Dual $\Z/4\Z$-code}

A $4$-frame in a lattice of rank $n$ is a collection of $n$
mutually orthogonal vectors of squared norm four.
The existence of such a frame realizes the lattice through
Construction A over $\Z/4\Z$.

\begin{proposition}\label{prop:frame}
The lattice $\Lambda$ contains a $4$-frame
$f_1,\ldots,f_{37}$ with $f_1=v_{00}$.
\end{proposition}
\begin{proof}
A computer search produces a $4$-frame, specified by nine
orbit representatives in the appendix.
Their orbits under $\langle C,M^6\rangle$ give $37$ mutually
orthogonal vectors of squared norm four.
Each lies in $V$ and pairs integrally with
$v_{ij},\kappa,\theta,\tau$.
Since these generate $\Lambda$ and $\Lambda$ is unimodular,
all $37$ vectors belong to $\Lambda$.
\end{proof}

The frame determines a code
\[
\mathcal C=
\left\{
\bigl((x,f_i)\bmod4\bigr)_{i=1}^{37}
:x\in\Lambda
\right\}
\subseteq(\Z/4\Z)^{37}.
\]
Recall that the Euclidean weight of a word is the sum of
the coordinate weights $0,1,4,1$ assigned to $0,1,2,3$,
respectively.

\begin{corollary}\label{cor:code}
The code $\mathcal C$ is a Type I self-dual
$\Z/4\Z$-code of length $37$ and minimum Euclidean
weight $16$. Moreover, its Construction A lattice is
isometric to $\Lambda$.
\end{corollary}
\begin{proof}
The map defining $\mathcal C$ has kernel
$\bigoplus_i\Z f_i$, whose index in the unimodular
lattice $\Lambda$ is $2^{37}$. Thus
$|\mathcal C|=2^{37}$.

For $x,y\in\Lambda$, orthogonality of the frame gives
\[
\sum_i(x,f_i)(y,f_i)=4(x,y)\equiv0\pmod4.
\]
Hence $\mathcal C$ is self-orthogonal, and its size
shows that it is self-dual.

Since $f_i/2$ form an orthonormal basis, these coordinates
identify $\Lambda$ with
\[
A_4(\mathcal C)
=\frac12\left\{
a\in\Z^{37}:a\bmod4\in\mathcal C
\right\}.
\]
The minimum norm four of $\Lambda$ therefore gives
$d_E(\mathcal C)\ge16$.
On the other hand, $v_{01}$ has squared norm four and
\[
(v_{01},f_1)=(v_{01},v_{00})=2,
\]
so it determines a nonzero codeword of Euclidean weight
at most $4v_{01}^2=16$.
Thus $d_E(\mathcal C)=16$.
Finally, $\mathcal C$ is Type I because $\Lambda$ is odd.
\end{proof}

\appendix
\section{Coordinates of a $4$-frame}\label{app:frame}

The frame in \Cref{prop:frame} is specified by nine symmetry
representatives. Write an ambient vector in block coordinates
as $(a;r,c,s)$, and let $r^\vee=(r_{-i})_i$ denote index
reversal, with the same notation for the other blocks.
The automorphisms $C$ and $M^6$ act by
\[
\begin{aligned}
C(a;r,c,s)&=(a;c,s^\vee,r^\vee),\\
M^6(a;r,c,s)&=(a;r^\vee,c^\vee,s^\vee).
\end{aligned}
\]
Applying $C^jM^{6\varepsilon}$, for $0\le j\le2$ and
$\varepsilon\in\{0,1\}$, to the nine vectors in
\Cref{tab:frame}, after dividing by $26$, gives
$37$ mutually orthogonal lines of squared norm four.
Their orbit sizes are
\[
1,\ 3,\ 3,\ 3,\ 3,\ 6,\ 6,\ 6,\ 6.
\]
Choosing one vector from each line gives the required
$4$-frame, with $f_1=v_{00}$.

\begin{table}[htbp]
\centering
\scriptsize
\setlength{\tabcolsep}{2.5pt}
\renewcommand{\arraystretch}{1.1}
\caption{Nine representatives of the $4$-frame,
multiplied by $26$. The column $a$ gives the
$e$-coordinate, and each block is indexed from $0$ to $12$.}
\label{tab:frame}
\begin{tabular}{r r c|rrrrrrrrrrrrr}
\hline
$k$ & $a$ & $X$
& 0 & 1 & 2 & 3 & 4 & 5 & 6 & 7 & 8 & 9 & 10 & 11 & 12\\
\hline
1 & 26 & $R$
& 26&0&0&0&0&0&0&0&0&0&0&0&0\\
  &    & $C$
& 26&0&0&0&0&0&0&0&0&0&0&0&0\\
  &    & $S$
& 26&0&0&0&0&0&0&0&0&0&0&0&0\\
\hline
2 & 13 & $R$
& -23&3&3&3&3&3&3&3&3&3&3&3&3\\
  &    & $C$
& 9&-17&9&-17&9&9&9&9&9&9&-17&9&-17\\
  &    & $S$
& 1&1&1&1&1&1&1&1&1&1&1&1&1\\
\hline
3 & 0 & $R$
& 0&0&0&0&0&26&-26&-26&26&0&0&0&0\\
  &   & $C$
& 0&0&0&0&0&0&0&0&0&0&0&0&0\\
  &   & $S$
& 0&0&0&0&0&0&0&0&0&0&0&0&0\\
\hline
4 & 0 & $R$
& 0&2&4&6&8&10&-14&14&-10&-8&-6&-4&-2\\
  &   & $C$
& 0&2&4&6&8&10&-14&14&-10&-8&-6&-4&-2\\
  &   & $S$
& 0&-2&-4&-6&-8&16&-12&12&-16&8&6&4&2\\
\hline
5 & 0 & $R$
& 0&4&8&12&-10&-6&-2&2&6&10&-12&-8&-4\\
  &   & $C$
& 0&4&8&-14&16&-6&-2&2&6&-16&14&-8&-4\\
  &   & $S$
& 0&-4&-8&14&10&6&2&-2&-6&-10&-14&8&4\\
\hline
6 & 0 & $R$
& 0&26&0&-26&0&0&0&0&0&0&0&0&0\\
  &   & $C$
& 0&0&-26&0&26&0&0&0&0&0&0&0&0\\
  &   & $S$
& 0&0&0&0&0&0&0&0&0&0&0&0&0\\
\hline
7 & 0 & $R$
& -14&-10&20&-2&2&6&10&-12&-8&-4&0&4&8\\
  &   & $C$
& 10&14&-8&-4&0&4&8&-14&-10&-6&-2&2&6\\
  &   & $S$
& 4&0&-4&-8&-12&10&6&2&-2&-6&-10&12&8\\
\hline
8 & 0 & $R$
& 6&-10&0&10&-6&4&14&-2&8&-8&2&-14&-4\\
  &   & $C$
& -8&2&-14&-4&6&-10&0&10&20&4&-12&-2&8\\
  &   & $S$
& 2&-8&8&-2&-12&4&-6&10&0&-10&6&-4&12\\
\hline
9 & 13 & $R$
& -11&-9&-7&-5&-3&-1&1&3&5&7&9&11&13\\
  &    & $C$
& -7&-5&-3&-1&1&3&5&7&9&-15&13&-11&17\\
  &    & $S$
& 5&3&1&-1&-3&-5&-7&-9&-11&13&11&9&7\\
\hline
\end{tabular}
\end{table}

\FloatBarrier
\section*{Acknowledgments}
I thank Ian Teixeira for introducing me to lattice theory and
explaining its connections with coding theory, which sparked
my interest in the subject.

\bibliographystyle{alpha}
\bibliography{biblio3}

\end{document}